\documentclass[12pt]{amsart}
\usepackage{amscd,amssymb}
\usepackage[graph,frame,poly,arc]{xy}  
\usepackage[plainpages,backref,urlcolor=blue]{hyperref}
\usepackage{color}
\usepackage{ textcomp }

\theoremstyle{plain}
\newtheorem{thm}[subsection]{Theorem}

\newtheorem{prop}[subsection]{Proposition}
\newtheorem{cor}[subsection]{Corollary}

\theoremstyle{definition}
\newtheorem{rk}[subsection]{Remark}
\newtheorem{definition}[subsection]{Definition}
\newtheorem{ex}[subsection]{Example}
\newtheorem{conj}[subsection]{Conjecture}

\numberwithin{equation}{section}
\newcommand{\A}{{\mathcal A}}

\newcommand{\al}{{\alpha}}
\newcommand{\be}{{\beta}}

\newcommand{\Q}{\mathbb{Q}}

\newcommand{\C}{\mathbb{C}}
\newcommand{\PP}{\mathbb{P}}

\newcommand{\N}{\mathbb{N}}

\DeclareMathOperator{\dd}{d}

\begin{document}

\title [Monodromy and pole order filtration on Milnor fiber]
{Computing the monodromy and pole order filtration on Milnor fiber cohomology of plane curves }

\author[Alexandru Dimca]{Alexandru Dimca}
\address{Universit\'e C\^ ote d'Azur, CNRS, LJAD, France }
\email{dimca@unice.fr}

\author[Gabriel Sticlaru]{Gabriel Sticlaru}
\address{Faculty of Mathematics and Informatics,
Ovidius University,
Bd. Mamaia 124, 
900527 Constanta,
Romania}
\email{gabrielsticlaru@yahoo.com }

\subjclass[2010]{Primary 32S40; Secondary 32S22, 32S55.}

\keywords{plane curve, Milnor fiber, monodromy, pole order filtration, b-function}

\begin{abstract} We describe an algorithm computing the monodromy and the pole order filtration on the Milnor fiber cohomology of any reduced projective plane curve $C$. The relation to the zero set of Bernstein-Sato polynomial of the defining homogeneous polynomial for $C$ is also discussed. When $C$ has some non weighted homogeneous singularities, then we have to assume that a conjecture holds in order to get some of our results. In all the examples computed so far this conjecture holds.

\end{abstract}
 
\maketitle


\section{Introduction} \label{sec:intro}

Let $C:f=0$ be a reduced plane curve  of degree $d\geq 3$  in the complex projective plane $\PP^{2}$, defined by a homogeneous polynomial $f \in S=\C[x,y,z]$.
Consider the corresponding complement $U=\PP^{2}\setminus C$, and the global Milnor fiber $F$ defined by $f(x,y,z)=1$ in $\C^3$ with monodromy action $h:F \to F$, 
$$h(x,y,z)=\exp(2\pi i/d)\cdot (x,y,z).$$ 
 To determine  the eigenvalues of the monodromy operators
\begin{equation} 
\label{mono1}
h^m: H^m(F,\C) \to H^m(F,\C)
\end{equation} 
for $m=1,2$ starting from $C$ or  $f$  is a rather difficult problem, going back to O. Zariski and attracting an extensive literature, see for instance \cite{A1, AD, BDS, CS, DHA, HE, L1, L2, OkaS, Deg, S2, D1}. When the curve $C:f=0$ is either free or nearly free, we have presented in \cite{DStFor} an efficient algorithm for listing the eigenvalues of the monodromy operator $h^1$, which determines completely the corresponding Alexander polynomial $\Delta_C(t)$, see Remark \ref{rkAlex} below for its definition.

In this paper we explain an approach  working in the general case. 
This time the results of our computation give not only the dimensions of the eigenspaces $H^m(F,\C)_{\lambda}$ of the monodromy, but also the dimensions of the graded
pieces $Gr^p_PH^m(F,\C)_{\lambda}$, where $P$ denotes the pole order filtration on $H^m(F,\C)$, see section 2 below for the definition. More precisely, the algorithm described here gives the following.
\begin{enumerate} 
\item the dimensions of the eigenspaces $H^m(F,\C)_{\lambda}$ for $m=1,2$  for any reduced curve  $C:f=0$, see Remark \ref{rkAlex}.
\item the dimensions of the graded
pieces $Gr^p_PH^1(F,\C)_{\lambda}$, for any reduced curve  $C:f=0$. Note that the $P^p$ filtration coincides to the Hodge filtration $F^p$ on $H^1(F,\C)$, see \cite[Proposition 2.2]{DStproj}.

\item the dimensions of the graded
pieces $Gr^p_PH^2(F,\C)_{\lambda}$, for a reduced curve  $C:f=0$ having only weighted homogeneous singularities. To achieve this efficiently one has to use the recent result by M. Saito stated below in Theorem \ref{thmconj0}. 

\item the dimensions of the graded
pieces $Gr^p_PH^2(F,\C)_{\lambda}$, for any reduced curve  $C:f=0$ under the assumption that a basic fact, stated as Conjecture \ref{conj30}, holds. This conjecture holds in all the examples we have computed so far, see Remark \ref{rkcheck}.

\end{enumerate}

The new information on the pole order filtration  $P$ can be applied to describe the set of roots of $b_f(-s)$, where $b_f(s)$ is the Bernstein-Sato polynomial of $f$, see for details \cite{Sa1}, \cite{Sa2}. In fact, using \cite[Theorem 2]{Sa1}, this comes down  to checking whether
$ Gr^p_PH^2(F,\C)_{\lambda}\ne 0,$ see Theorem \ref{thmBS} for a precise statement and our applications described in Corollaries \ref{corAA'},  \ref{corAB'},  \ref{corAC'}.

Here is in short how we proceed. Let $\Omega^j$ denote the graded $S$-module of (polynomial) differential $j$-forms on $\C^3$, for $0 \leq j \leq 3$. The complex $K^*_f=(\Omega^*, \dd f \wedge)$ is just the Koszul complex in $S$ of the partial derivatives $f_x$, $f_y$ and $f_z$ of the polynomial $f$. The general theory says that there is a spectral sequence $E_*(f)$, whose first term $E_1(f)$ is computable from the cohomology of the Koszul complex $K^*_f$ and whose limit
$E_{\infty}(f)$ gives us the action of monodromy operator  on the graded pieces of the reduced cohomology $\tilde H^*(F,\C)$ of $F$ with respect to the pole order filtration $P$, see \cite{Dcomp}, \cite[Chapter 6]{D1}, \cite{DS1}.

Our approach takes a simpler form when $C$ is assumed to have only weighted homogeneous singularities, e.g.  when $C$ is a line arrangement $\A$.
This comes from the following result due to M. Saito \cite{Sa3}, see for a more precise statement Theorem \ref{thmconj} below.

\begin{thm}
\label{thmconj0}
If the reduced plane curve $C:f=0$ has only weighted homogeneous singularities, then
the spectral sequence $E_*(f)$ degenerates at the $E_2$-term. 
\end{thm}
This result has been conjectured already in \cite{Dcomp} and has been checked in many cases using a computer program in \cite{DStFor}. The converse implication is known to hold, see \cite[Theorem 5.2]{DS1}.  
The algorithm described in \cite{DStFor} for free and nearly free curves actually computes the $E_2$-term of this spectral sequence.

In this paper we modify the algorithm in \cite{DStFor} such that it applies to any reduced curve. First we compute (a large part) of the $E_2-$term in Section 4, which is (more than) enough when $C$ has only weighted homogeneous singularities, see Remark \ref{rkwh}, or when we are interested only in the Alexander polynomial $\Delta_C(t)$.
Then in Section 5  we compute  the relevant part of the $E_3$-term of the above spectral sequence. Conjecture \ref{conj30} tells us that essentially $E_3=E_{\infty}$ and that we have missed no information on the $P$-filtration on $H^2(F,\C)$ by looking only at some of the terms in the $E_2$ and $E_3$ pages.

The computations in this note were made using the computer algebra system  Singular \cite{Sing}. 
The corresponding codes are available at \\
\url{http://math.unice.fr/~dimca/singular.html}
\medskip

We thank Morihiko Saito for his remarks that greatly improved the clarity of the presentation. In particular, Conjecture \ref{conj30} and Proposition \ref{propMS} were suggested by him, and they provide an efficient method to certify all the results given by our algorithm. For an alternative approach, see Morihiko Saito's paper \cite{Sa4}.

\section{Gauss-Manin complexes, Koszul complexes, and Milnor fiber cohomology} \label{sec2}

Let $S$ be the polynomial ring $\C[x,y,z]$ with the usual grading and consider a reduced homogeneous  polynomial $f \in S$ of degree $d$. The graded Gauss-Manin complex $C_f^*$ associated to $f$ is defined by taking $C_f^j=\Omega^j[\partial_t]$, i.e. formal polynomials in $\partial_t$ with coefficients in the space of differential forms $\Omega^j$, where $\deg \partial_t=-d$ and the differential
$\dd: C_f^j \to C_f^{j+1}$ is $\C$-linear and given by 
\begin{equation} 
\label{difC}
 \dd (\omega \partial_t^q)=(\dd \omega)\partial_t^q-(\dd f \wedge \omega) \partial_t^{q+1},
\end{equation} 
see for more details \cite[Section 4]{DS1}. 
The complex $C_f^*$ has a natural increasing filtration $P'_*$ defined by
\begin{equation} 
\label{filC}
 P'_qC^j_f=\oplus_{i \leq q+j}\Omega^j\partial_t^i.
\end{equation} 
If we set $P'^q=P'_{-q}$ in order to get a decreasing filtration, then one has
\begin{equation} 
\label{grC}
 Gr^q_{P'}C^*_f=\sigma_{\geq q}(K^*_f((3-q)d)),
\end{equation} 
the truncation of a shifted version of the Koszul complex $K^*_f$, where $\sigma$ denotes the stupid filtration, see for instance \cite[Remark 1.1.15]{D2}.
Moreover, this yields a decreasing filtration $P'$ on the cohomology groups $H^j(C^*_f)$ and a spectral sequence
\begin{equation} 
\label{spsqC}
 E_1^{q,j-q}(f) \Rightarrow H^j(C^*_f).
\end{equation} 
On the other hand, the reduced cohomology $\tilde H^j(F,\C)$ of the Milnor fiber $F:f(x,y,z)=1$ associated to $f$ has a pole order decreasing filtration $P$, see \cite[Section 3]{DS1}, such that there is a natural identification for any integers $q$, $j$ and $k \in [1,d]$
\begin{equation} 
\label{filH}
 P'^{q+1}H^{j+1}(C^*_f)_k=P^q\tilde H^j(F,\C)_{\lambda},
\end{equation} 
where $\lambda=\exp (-2 \pi ik/d).$ 
Since the Milnor fiber $F$ is a smooth affine variety, its cohomology groups $H^m(F, \C)$ have a decreasing Hodge filtration $F$ coming from the mixed Hodge structure constructed by Deligne, see \cite{PeSt}. The two filtrations $P$ and $F$ are related by the inclusion
\begin{equation} 
\label{PFincl}
 F^sH^{m}(F,\C) \subset P^s  H^m(F,\C),
\end{equation} 
for any integers $s,m$, see formula $(4.4.8)$ in \cite{DS1}. Note also that 
\begin{equation} 
\label{PFincl1}
F^0H^{m}(F,\C)= H^m(F,\C) \text{ and } P^{m+1}  H^m(F,\C)=0, 
\end{equation}
for any integer $m$, where the first equality comes from the general properties of the Hodge filtration, see \cite{PeSt}, while the second equality follows from the definition of the filtration $P$.

The $E_1$-term of the spectral sequence 
\eqref{spsqC} is completely determined by the morphism of graded $\C$-vector spaces
\begin{equation} 
\label{diff1}
 \dd ' : H^2(K^*_f) \to H^3(K^*_f),
\end{equation} 
induced by the exterior differentiation of forms, i.e. $\dd ' :[\omega] \mapsto [\dd (\omega)]$.
Note that this morphism $\dd'$ coincides with the morphism $\dd ^{(1)}: N \to M$ considered in 
\cite{DS1}, up-to some shifts in gradings. More precisely, if we set for $k \in [1,d]$,
\begin{equation} 
\label{newspsq}
E_1^{s,t}(f)_k=H^{s+t+1}(K^*_f)_{td+k} \text{ and } \dd_1:E_1^{s,t}(f)_k \to E_1^{s+1,t}(f)_k, \dd _1:[\omega] \mapsto [\dd (\omega)],
\end{equation} 
we get a new form of (a homogeneous component of) the spectral sequence \eqref{spsqC},
where $E_1^{s,t}(f)_k=0$ for $s+t \notin \{1,2\}$.  Hence, with the above notation, one has
\begin{equation} 
\label{limit}
E_{\infty}^{s,t}(f)_k=Gr_P^s\tilde H^{s+t}(F,\C)_{\lambda }.
\end{equation} 
The case $k=d$ is also discussed in \cite{DStEdin}, see also \cite{Sa4}. The second differential 
\begin{equation} 
\label{diff2}
 \dd_2:E_2^{s,t}(f)_k \to E_2^{s+2,t-1}(f)_k, 
\end{equation}
for $s+t=1$ can be described as follows. If $[\omega] \in E_2^{s,t}(f)_k$, then $\dd_1([\omega])=0$, which means that $\dd f \wedge \omega=0$ and there is a 2-form $\al \in \Omega^2$ such that $\dd (\omega)=\dd f \wedge \al$. Then one has
$$ \dd_2([\omega])=[\dd \al].$$

One has the following key result due to M. Saito \cite{Sa3}.

\begin{thm}
\label{thmconj}
The spectral sequences \eqref{newspsq} degenerate at the $E_2$-term when the reduced plane curve $C:f=0$ has only weighted homogeneous singularities.
\end{thm}

We need the following result, see  \cite[Theorem 5.3]{DS1}. Let $\al_{p_i,j}$ with $j=1,...,\mu(C,p_i)$ be the spectral numbers of the plane curve singularity $(C,p_i)$, where each spectral number is repeated as many times as its multiplicity in the spectrum of $(C,p_i)$.

\begin{thm}
\label{thminj}
Assume that the plane curve $C:f=0$ in $\PP^2$ has only weighted homogeneous singularities and consider
$$E_2^{1-t,t}(f)_k=\ker \left \{ \dd ' : H^2(K^*_f)_{td+k} \to H^3(K^*_f)_{td+k}\right \}.$$
Then $\dim E_2^{1-t,t}(f)_k  \leq N(C,t,k)$, where $N(C,t,k)$ is the number of spectral numbers $\al_{p_i,j}$  equal to $\frac{td+k}{d}$,
when $p_i$ ranges over all the singularities of the curve $C$ and $j=1,...,\mu(C,p_i).$
\end{thm}

Theorem \ref{thminj} also implies the following weaker form of Theorem \ref{thmconj}, see also \cite[Corollary 5.5]{DS1}.

\begin{cor} 
\label{corvanishing} 
 If the degree $d$ reduced plane curve $C$ has only weighted homogeneous singularities, then 
$E_{2}^{1,0}(f)_k=0$ for $k<\alpha(C)d$, where $\al(C)$ is the minimum of the spectral numbers
of the singularities of $C$. Moreover, the spectral sequences  \eqref{newspsq} degenerate at the $E_3$-term and the only possibly non-zero differentials in the $E_2$-terms are the differentials
$$d_2: E_2^{0,1}(f)_k \to E_2^{2,0}(f)_k$$
for $k=1,...,(1-\alpha(C))d,$ i.e. one has $E_2^{1-t,t}(f)_k=0$ for $q=td+k > (2-\alpha(C))d$.
\end{cor}

\begin{cor}
\label{corAD}
Let $C: f=0$ be a reduced plane curve of degree $d$ having only weighted homogeneous singularities.
Then $\dim Gr_P^{1} H^{1}(F,\C)_{\lambda }= 0$ for $k <\al(C)d$ and $\dim Gr_P^{0} H^{1}(F,\C)_{\lambda }= 0$
for $k > (1-\al(C))d$.
\end{cor}

Corollary \ref{corvanishing} implies that the following three related conjectures hold when all the singularities of $C$ are weighted homogeneous. 

\begin{conj}
\label{conj3}
The spectral sequences  \eqref{newspsq} degenerates at the $E_3$-term for any reduced plane curve $C:f=0$.
\end{conj} 
We have no idea how to prove this conjecture, not even how to check that it holds on a specific example.
For our practical purposes, we need only the following weaker version of this conjecture.
This version can also be checked on a given example using the algorithm described below, see Proposition \ref{propMS}.
\begin{conj}
\label{conj30} In the spectral sequences  \eqref{newspsq}  of any reduced plane curve $C:f=0$,
the following equalities hold
$$\dim E_{3}^{2-t,t}(f)_k =\dim E_{\infty}^{2-t,t}(f)_k,$$
for any integer $k \in [1,d]$ and any $t=0,1,2$.

\end{conj} 
To speed up the computations, the following more precise form of the above claims would be very useful. 
\begin{conj}
\label{conjmu-tau}
For any reduced plane curve $C:f=0$, 
there is a positive  integer $q_0 \leq 3d+1$, such that
 $$\dim E_2^{1-t,t}(f)_{k}= \mu(C)-\tau(C) \text{ and }  E_3^{1-t,t}(f)_{k}=0$$
 for any $q=td+k\geq q_0$, where $\mu(C)$ is the total Milnor number of $C$, that is the sum of all
the  Milnor numbers of the singular points of $C$. Similarly, $\tau(C)$ is the total Tjurina number of $C$, that is the sum of all
the Tjurina numbers of the singular points of $C$.
The minimal positive integer $q_0$ satisfying the above property, if it exists, is denoted by $q_0(f)$. 
\end{conj}
 
It is clear that Conjecture \ref{conjmu-tau} implies Conjecture \ref{conj30}, but again
we have no idea how to prove Conjecture \ref{conjmu-tau}, not even how to check that it holds on a specific example.

\begin{rk}
\label{rkq0}
Corollary \ref{corvanishing} implies that when all the singularities of $C$ are weighted homogeneous, one has $q_0(f) \leq (2-\alpha(C))d+1$. Note that by the semicontinuity of the spectrum, one has $\alpha(C) \geq 2/d$ and hence 
the vanishing $E_2^{1-t,t}(f)_k=0$ holds in particular for $q=td+k \geq 2d-1$.

However, there are curves for which $q_0(f)>2d-1$,  see Example \ref{exhighq0}.

\end{rk}

\begin{rk}
\label{rkA}
The two filtrations $F$ and $P$ coincide on $H^{1}(F,\C)_{\lambda }$ always, see \cite[Proposition 2.2]{DStproj}.
  On the other hand, the two filtrations $F$ and $P$ do not coincide on $H^{2}(F,\C)_{\lambda }$
even in very simple cases, e.g. $C:f=(x^2-y^2)(x^2-z^2)(y^2-z^2)=0$ and ${\lambda }=-1$.
A computation of the Hodge filtration on $H^{2}(F,\C)$ in this case can be found in \cite{PBthesis, DHA}. Note  also that the mixed Hodge structure on 
$H^2(F,\C)_{\ne 1}$ is not pure in general. For a line arrangement, one can use the formulas for the spectrum given in \cite{BS} to study the interplay between monodromy and Hodge filtration on 
$H^2(F,\C)_{\ne 1}$, see \cite[Remark 2.5]{DStproj}.

\end{rk}

\section{Jacobian syzygies of plane curves}

Consider the graded $S-$submodule $AR(f) \subset S^{3}$ of {\it all relations} involving the derivatives of $f$, namely
$$\rho=(a,b,c) \in AR(f)_q$$
if and only if  $af_x+bf_y+cf_z=0$ and $a,b,c$ are in $S_q$, the space of homogeneous polynomials of degree $q$.   Let $d_1=mdr(f)$ be the minimal degree of a relation in $AR(f)$. We assume in the sequel that $d_1>0$, which is equivalent to saying that $C$ is not the union of $d$ lines passing through one point, a case easy to handle directly.

To each syzygy $\rho=(a,b,c) \in AR(f)_q$ we associate a differential $2$-form
\begin{equation} 
\label{form}
\omega(\rho)=a \dd y \wedge \dd z -b  \dd x \wedge \dd z+ c \dd x \wedge \dd y\in \Omega^2_{q+2}
\end{equation} 
such that the relation $af_x+bf_y+cf_z=0$ becomes $\dd f \wedge \omega(\rho)=0.$
A relation $\rho=(a,b,c)$ is called a Koszul relation if the associated differential form
$\omega (\rho)$ is of the form $\dd f \wedge \alpha$, for some $\alpha \in \Omega ^1$.
The set $KR(f)$ of these Koszul relations forms a graded $S$-submodule in $AR(f)$.

Hence, up to a shift in degrees, for any polynomial $f$ there is an identification 
$$AR(f)(-2)=Syz(f):=\ker \{ \dd f \wedge: \Omega ^2 \to \Omega ^3 \},$$
such that the Koszul relations $KR(f)$ inside $AR(f)$ correspond to the submodule $\dd f \wedge \Omega^1$ in $Syz(f)$. Since $C:f=0$ has only isolated singularities, it follows that $H^1(K^*_f)=0$, i.e.
the following sequence, where the morphisms are the wedge product by $\dd f$, is exact for any $j$
$$ 0 \to \Omega^0_{j-2d} \to \Omega^1_{j-d} \to (\dd f \wedge \Omega^1)_j \to 0.$$
In particular, one has
\begin{equation} 
\label{dimKR1}
\dim (\dd f \wedge \Omega^1)_j =0 \text{  for } j \leq d,
\end{equation} 

\begin{equation} 
\label{dimKR2}
\dim (\dd f \wedge \Omega^1)_j =3 {j-d+1 \choose 2} \text{  for } d <j < 2d.
\end{equation} 
and
\begin{equation} 
\label{dimKR3}
\dim (\dd f \wedge \Omega^1)_j =3 {j-d+1 \choose 2} -{j-2d+2 \choose 2} \text{  for } j  \geq 2d.
\end{equation}

\begin{rk}
\label{rkdim}
Let $J_f$ be the Jacobian ideal spanned by $f_x,f_y,f_z$ in $S$, and denote by $M(f)=S/J_f$ the corresponding Jacobian (or Milnor) algebra of $f$. Let $m(f)_j= \dim M(f)_j$ for $j \geq 0$ and recall the formulas
\begin{equation} 
\label{dimER1}
\dim H^2(K^*_f)_j = m(f)_{j+d-3}-m(f_s)_{j+d-3} \text{  for } 2 \leq j \leq 2d-3,
\end{equation} 
and $\dim H^2(K^*_f)_j =\tau(C)$ for $j \geq 2d-2$, where $C_s:f_s=0$ denotes a smooth curve of degree $d$, see \cite{DBull}.
Since $H^2(K^*_f)_j = Syz(f)_j/(\dd f \wedge \Omega^1)_j$, the combination of the formulas \eqref{dimKR1}, \eqref{dimKR2}, \eqref{dimKR3} and \eqref{dimER1} above gives us formulas for the dimensions 
$syz(f)_j= \dim Syz(f)_j$ for any $j\geq 3$. Note that $Syz(f)_j=AR(f)_{j-2}=0$ for any $j< 3$ by our assumption $d_1=mdr(f)>0$.
\end{rk}
\begin{definition}
\label{def}
For a reduced plane curve $C:f =0$ of degree $d$, we recall the following invariants.

\noindent (i) the {\it coincidence threshold} 
$$ct(f)=\max \{q:m(f)_k=m(f_s)_k \text{ for all } k \leq q\},$$
with $f_s$  a homogeneous polynomial in $S$ of the same degree $d$ as $f$ and such that $C_s:f_s=0$ is a smooth curve in $\PP^2$.

\noindent (ii) the {\it stability threshold} 
$st(f)=\min \{q~~:~~m(f)_k=\tau(C) \text{ for all } k \geq q\}.$

\end{definition}

\section{The first cycle in the algorithm: computing $E^{s,t}_2(f)_k$ }

\subsection{The computation of  $E^{s,t}_2(f)_k$ for $s+t=1$, $0 \leq t \leq 3$}
For $3 \leq q \leq 4d$, we set $q_1=q-d$ and consider the linear mapping
\begin{equation} 
\label{fi4}
\phi'_q:S^3_{q-2} \times S^3_{q_1-2}\to S_{q-3+d} \times S_{q-3},
\end{equation} 
given by 
$$((a,b,c),(u,v,w)) \mapsto (af_x+bf_y+cf_z, a_x+b_y+c_z-uf_x-vf_y-wf_z).$$ 
Since $S_j=0$ for $j<0$, this map has a simpler form for $q<d+2$. This map puts together the differentials $\dd_1$ and $\dd_2$ from the spectral sequence \eqref{newspsq}. Indeed, if $\omega \in \Omega^2$ (resp. $\al \in \Omega^2$) is given by the formula \eqref{form} starting with the triple  $(a,b,c) \in S^3_{q-2}$ (resp. the triple $(u,v,w) \in S^3_{q_1-2}$), one sees that essentially $\phi'_q((a,b,c),(u,v,w))$ corresponds to the pair
$$(\dd f \wedge \omega, \dd \omega-\dd f \wedge \al).$$ 
It is clear that 
$$((a,b,c),(u,v,w)) \in K'_q:=\ker \phi'_q$$
if and only if $\dd f \wedge \omega=0$ and $ \dd \omega= \dd f \wedge \al$. Note that if $\omega=\dd f \wedge \eta$, one can take $\al=-\dd \eta$ and the corresponding pair $(\omega,\al)$ gives rise to an element in $K'_q$.
Consider the projection $B_q \subset S^3_{q-2}$ of $K'_q$ on the first component
and note that $B_q/(\dd f \wedge \Omega^1)_q$ can be identified to $E_2^{1,0}(f)_{q}$, for $3 \leq q \leq d$, to $E_2^{0,1}(f)_{q-d}$, for $d+1 \leq q \leq 2d$, to $E_2^{-1,2}(f)_{q-2d}$, for $2d+1 \leq q \leq 3d$,
and respectively to $E_2^{-2,3}(f)_{q-3d}$, for $3d+1 \leq q \leq 4d$. 

On the other hand, the kernel of the projection $K'_q \to B_q$ can be identified to the set of forms $\al' \in \Omega^2$ such that $\dd f \wedge \al'=0$. It follows that if we set $k'_q= \dim K'_q$ and
$$\epsilon'_q= k'_q -syz(f)_{q_1} -\dim (\dd f \wedge \Omega^1)_q,$$
we have
$\epsilon'_q=\dim E_2^{1,0}(f)_{q}$ for $3 \leq q \leq d$,  $\epsilon'_q=\dim E_2^{0,1}(f)_{q-d}$ for $d+1 \leq q \leq 2d$,  $\epsilon'_q=\dim E_2^{-1,2}(f)_{q-2d}$ for $2d+1 \leq q \leq 3d$, and respectively $\epsilon'_q=\dim E_2^{-2,3}(f)_{q-3d}$ for $3d+1 \leq q \leq 4d$

By convention, we set $K'_j=0$ and $k'_j=\epsilon'_j=0$ for $j=0,1,2$.

\subsection{The computation of  $E^{s,t}_2(f)_k$ for $s+t=2$, $0 \leq t \leq 3$}
Recall that, setting $q=td+k$, we have $1 \leq q \leq 4d $ and
$$\dim E^{2-t,t}_1(f)_k= \dim H^3(K^*_f)_{q} = \dim M(f)_{q-3}= m(f)_{q-3}.$$
It follows that
$$\theta_q:=\dim E^{2-t,t}_2(f)_k=\dim E^{2-t,t}_1(f)_k-(\dim E^{1-t,t}_1(f)_k-\dim E^{1-t,t}_2(f)_k)=$$
$$=m(f)_{q-3}-syz(f)_q+\dim (\dd f \wedge \Omega^1)_q+\epsilon'_q.$$
This can be rewritten as
$$\theta_q=m(f)_{q-3}-m(f)_{q+d-3}+m(f_s)_{q+d-3}+\epsilon'_q $$  for $2 \leq q \leq 2d-3$
and
$$\theta_q=m(f)_{q-3}-\tau(C)+\epsilon'_q $$  for $2d-2 \leq  q \leq 3d-4$, resp. $\theta_q=\epsilon'_q $ for $3d-3 \leq q \leq 4d$
in view of Remark \ref{rkdim} above and since $st(f) \leq 3d-6$.
\begin{rk}
\label{rkwh}
When the degree $d$ reduced plane curve $C$ has only weighted homogeneous singularities, we have
$$\epsilon'_q =0$$
for any $q > (2-\al(C))d$ as implied by Corollary \ref{corvanishing}. Moreover, Theorem \ref{thmconj0} implies that in this case 
$$E^{s,t}_2(f)_k=E^{s,t}_{\infty}(f)_k$$
and hence to determine the monodromy action and the pole order filtration on Milnor fiber cohomology in this case it is enough to compute the integers $\epsilon'_q$ only for $q \leq 2d-2$. Indeed, all the terms $E^{s,t}_2(f)_k=E^{s,t}_{\infty}(f)_k$ for $s+t=1$ and $q>2d-2$  are trivial, as it follows from \eqref{limit} and the properties of the pole order filtration $P$, see  \eqref{PFincl} and \eqref{PFincl1}.
Therefore for a reduced plane curve $C$ having only weighted homogeneous singularities the algorithm stops at this stage. Note that when $C$ is a line arrangement one can compute only  $\theta_q$ for $q \leq 2d$, since $P^1H^2(F,\C)=H^2(F,\C)$ as shown in Corollary \ref{corAA}, but in general one may need some values $\theta_q$ for $2d< q \leq st(f)+2$ as shown by Corollary \ref{corAC} and Example \ref{extc2,4}, where $P^1H^2(F,\C) \ne H^2(F,\C)$.

\end{rk}

\begin{rk}
\label{rkAlex}

One can consider  the characteristic polynomials of the monodromy, namely
\begin{equation} 
\label{Delta}
\Delta^j_C(t)=\det (t\cdot Id -h^j|H^j(F,\C)),
\end{equation} 
for $j=0,1,2$. It is clear that, when the curve $C$ is reduced,  one has $\Delta^0_C(t)=t-1$, and moreover
\begin{equation} 
\label{Euler}
\Delta^0_C(t)\Delta^1_C(t)^{-1}\Delta^2_C(t)=(t^d-1)^{\chi(U)},
\end{equation} 
where $\chi(U)$ denotes the Euler characteristic of the complement $U$, see for instance \cite{L1, LV} or \cite[Proposition 4.1.21]{D1}. Since
 $$\chi(U)=(d-1)(d-2)+1-\mu(C)$$
it follows that the polynomial $\Delta_C(t)=\Delta^1_C(t)$, also called the Alexander polynomial of $C$, see \cite{R}, determines the remaining polynomial $\Delta^2_C(t)$. Note that the computation of the dimension of the term $E^{1,0}_2(f)_k$, for $k=1,...,d$  described in this section  is enough to determines the Alexander polynomial $\Delta_C(t)$ for all reduced plane curves. Indeed, the equality of the pole order filtration $P^p$ with  the Hodge filtration $F^p$ on $H^1(F,\C)$, see \cite[Proposition 2.2]{DStproj} and the obvious equality 
$$\dim H^{0,1}(F,C)_{\lambda}=\dim H^{1,0}(F,C)_{\bar \lambda},$$
with $ \bar \lambda$ the complex conjugate of $ \lambda$, imply that
$$\dim H^{1,0}(F,C)_{\lambda}= \dim E^{1,0}_2(f)_k \text{ and }\dim H^{0,1}(F,C)_{\lambda}= \dim E^{1,0}_2(f)_{d-k},$$
where $\lambda=\exp(-2\pi ik/d)$.

Note that the results by H. Esnault in \cite{HE} and even more clearly the results of E. Artal Bartolo in \cite{A1} can be used to compute the Alexander polynomials in many cases. However these results do not seem to be easily implementable as algorithms to performed computer aided computations.

We also note that our results, reducing the computation of the Alexander polynomial to linear algebra, imply in a clear way that this polynomial is not changed when we apply an automorphism of $\C$ over $\Q$ to our defining polynomial $f$. In particular, the Alexander polynomial cannot be used to distinguished Galois conjugate Zariski pairs as noticed already in \cite{A2}.
\end{rk}

\begin{rk}
\label{rkmu-tau}
In practice, the computation of $\epsilon'_q=\dim E_2^{-1,2}(f)_{q-2d}$ for $2d+1 \leq q \leq 4d$
takes a lot of time. A way to estimate in general the numerical invariant $q_0(f)$ introduced in Conjecture \ref{conjmu-tau}  would be of a great help.
All the examples computed so far suggest that  Conjecture \ref{conjmu-tau} holds. Some of these examples are given in the final section.
\end{rk}

\section{The second cycle in the algorithm: computing $E^{s,t}_3(f)_k$}
In this section we explain how the algorithm described in the previous section has to be continued in the case of the presence of non weighted homogeneous singularities, and if one wants to control the $P$-filtration on $H^2(F,\C)$.
The general construction is described in the subsection \ref{nwh4},
but for the clarity of exposition and the optimization of the computer time we discuss 
several cases.

\subsection{The case $3 \leq q \leq d+1$} \label{nwh1}
We consider the linear mapping
\begin{equation} 
\label{fi1}
\phi_q:S^3_{q-2} \to S_{q-3+d} \times S_{q-3},
\end{equation} 
given by 
$$(a,b,c) \mapsto (af_x+bf_y+cf_z, a_x+b_y+c_z).$$ 
This map is a simpler version of the map in \eqref{fi4}, and in fact $\phi_q=\phi'_q$ in this range.

It follows that $(a,b,c)$ is in the kernel $K_q:= \ker \phi_q$ if and only if the corresponding form $\omega$ is in $E_2^{1,0}(f)_q$ for $q=3,...,d$ or in $E_2^{0,1}(f)_1$ for $q=d+1$.
Note that one obviously has $E_2^{1,0}(f)_q=E_{\infty}^{1,0}(f)_q $ for $q=3,...,d$, and also
$E_2^{0,1}(f)_1=E_{\infty}^{0,1}(f)_1$ as explained in \cite[Remark 2.4]{DStFor}.
We set $k_q=\dim K_q$ and $\epsilon_q= k_q - \dim (\dd f \wedge \Omega^1)_q$. Since the class $[\omega]$ of the form $\omega$ is determined up to an element of $(\dd f \wedge \Omega^1)_q$, it follows that
$$\epsilon_q=\dim E_{\infty}^{1,0}(f)_q = \dim Gr_P^1 H^{1}(F,\C)_{\lambda }$$ 
for $q=3,...,d$ and $\epsilon_{d+1}=\dim E_{\infty}^{0,1}(f)_q= \dim Gr_P^0 H^{1}(F,\C)_{\lambda } $, with $\lambda = \exp(-2\pi iq/d)$. By convention, we set $K_j=0$ and $k_j=\epsilon_j=0$ for $j=0,1,2$.

\subsection{The case $d+2 \leq q \leq d+d_1+1$ where $d_1=mdr(f)$}  \label{nwh2}
We set $q_1=q-d$ and 
we consider the linear mapping
\begin{equation} 
\label{fi2}
\phi_q:S^3_{q-2} \times S^3_{q_1-2}\to S_{q-3+d} \times S_{q-3}\times S_{q_1-3},
\end{equation} 
given by 
$$((a,b,c),(u,v,w)) \mapsto (af_x+bf_y+cf_z, a_x+b_y+c_z-uf_x-vf_y-wf_z, u_x+v_y+w_z).$$ 
This map can be regarded as an extension of the map \eqref{fi4} obtained by adding the last component.
 Then it is clear that 
$$((a,b,c),(u,v,w)) \in K_q=\ker \phi_q$$
if and only if $\dd f \wedge \omega=0$, $ \dd \omega= \dd f \wedge \al$ and $\dd \al=0$. Since $q_1 \leq d_1+1$, it follows that $E_1^{1,0}(f)_{q_1}=0$, and hence $E_2^{2,0}(f)_{q_1}=E_1^{2,0}(f)_{q_1}=M(f)_{q_1-3}$. Consider the projection $A_q \subset S^3_{q-2}$ of $K_q$ on the first component
and note that $A_q/(\dd f \wedge \Omega^1)_q$ can be identified to $E_3^{0,1}(f)_{q_1}$. Moreover
$E_3^{0,1}(f)_{q_1}=E_{\infty}^{0,1}(f)_{q_1}$ since $q_1 \leq d$.
The kernel of the projection $K_q \to A_q$ can be identified to the set of forms $\al' \in \Omega^2$ such that $\dd f \wedge \al'=0$ and $\dd \al'=0$. Since $q_1\leq d_1+1$, the fist condition implies $\al'=0$, and hence $k_q=\dim K_q=\dim A_q$ in this case.
It follows that if we set
$$\epsilon_q= k_q -\dim (\dd f \wedge \Omega^1)_q,$$
we have again
$$\epsilon_q=\dim E_{\infty}^{0,1}(f)_{q _1}= \dim Gr_P^0 H^{1}(F,\C)_{\lambda }.$$
\subsection{The case $d+d_1+2 \leq q \leq 2d $ where $d_1=mdr(f)$}  \label{nwh3}
We set again $q_1=q-d$ and 
we consider the linear mapping
\begin{equation} 
\label{fi3}
\phi_q:S^3_{q-2} \times S^3_{q_1-2}  \times S^3_{q_1-2} \to S_{q-3+d} \times S^2_{q-3}\times S_{q_1-3},
\end{equation} 
given by 
$$((a,b,c),(u,v,w), (u',v',w')) \mapsto (\Phi_1,\Phi_2,\Phi_3,\Phi_4),$$
where
$\Phi_1=af_x+bf_y+cf_z$, $\Phi_2= a_x+b_y+c_z-uf_x-vf_y-wf_z$, $\Phi_3=u'f_x+v'f_y+w'f_z$
and $\Phi_4=
u_x+v_y+w_z-u'_x-v'_y-w'_z.$
With the  notation from the previous case, let $\be \in \Omega^2$ be the form associated to the triple $(u',v',w')$. Then it is clear that
$$((a,b,c),(u,v,w),(u',v',w')) \in K_q=\ker \phi_q$$
if and only if $\dd f \wedge \omega=0$, $ \dd \omega= \dd f \wedge \al$, $\dd f \wedge \be=0$ and $\dd \al=\dd \be$. Since $q_1 \leq d$, it follows that $\Phi_4=0$ has the same meaning in both $S_{q_1-3}$ and in $M(f)_{q_1-3}=S_{q_1-3}.$

 Consider the projection $A_q \subset S^3_{q-2}$ of $K_q$ on the first component
and note that $A_q/(\dd f \wedge \Omega^1)_q$ can be identified to $E_3^{0,1}(f)_{q_1}$, which again, is clearly
the same as $E_{\infty}^{0,1}(f)_{q_1}.$ The kernel of the projection $K_q \to A_q$ can be identified to the set of forms $\al', \be' \in \Omega^2$ such that $\dd f \wedge \al'=\dd f \wedge \be'=0$ and $\dd \al'=\dd \be'$. By setting $\gamma= \al'-\be'$, we see that this is the same as the set of forms $\al', \gamma \in \Omega^2$ such that $\dd f \wedge \al'=\dd f \wedge \gamma=0$ and $\dd \gamma=0$. This says exactly that $\al' \in Syz(f)_{q_1}$ and $\gamma \in K_{q_1}=K'_{q_1}$, where $K'_{q_1}$ was introduced in the first cycle above.
It follows that if we set
$$\epsilon_q= k_q -syz(f)_{q_1}-k'_{q_1} -\dim (\dd f \wedge \Omega^1)_q,$$
we have again
$$\epsilon_q=\dim E_{\infty}^{0,1}(f)_{q _1}= \dim Gr_P^0 H^{1}(F,\C)_{\lambda }.$$
Note that the value for  $syz(f)_{q_1}$ is determined in Remark \ref{rkdim}, the value for
$k_{q_1}$ is computed in the first step of our algorithm, and the value of $\dim (\dd f \wedge \Omega^1)_q$ is given in the equations \eqref{dimKR1}, \eqref{dimKR2}, and \eqref{dimKR3}.
\subsection{The case $2d+1 \leq q \leq 4d $ }  \label{nwh4}
We set  $q_1=q-d$ and $q_2=q-2d$, and
we consider the linear mapping
\begin{equation} 
\label{fi31}
\phi_q:S^3_{q-2} \times S^3_{q_1-2}  \times S^3_{q_1-2} \times S^3_{q_2-2} \to S_{q-3+d} \times S^2_{q-3}\times S_{q_1-3},
\end{equation} 
given by 
$$((a,b,c),(u,v,w), (u',v',w'), (u'',v'',w'')) \mapsto (\Phi_1,\Phi_2,\Phi_3,\Phi_4),$$
where
$\Phi_1=af_x+bf_y+cf_z$, $\Phi_2= a_x+b_y+c_z-uf_x-vf_y-wf_z$, $\Phi_3=u'f_x+v'f_y+w'f_z$
and $\Phi_4=
u_x+v_y+w_z-u'_x-v'_y-w'_z-u''f_x-v''f_y-w''f_z.$
With the  notation from the previous case, let $\eta \in \Omega^2$ be the form associated to the triple $(u'',v'',w'')$. Then it is clear that
$$((a,b,c),(u,v,w),(u',v',w'), (u'',v'',w'')  ) \in K_q=\ker \phi_q$$
if and only if $\dd f \wedge \omega=0$, $ \dd \omega= \dd f \wedge \al$, $\dd f \wedge \be=0$ and $\dd \al-\dd \be= \dd f \wedge \eta$. As above, one considers the projection
$K_q \to A_q$ and note that the kernel of this projection can be identified to the space $K'_{q_1}$, whose dimension $k'_{q_1}$ was determined in the first cycle above. Hence
$$\epsilon_q= k_q -syz(f)_{q_1}-k'_{q_1} -\dim (\dd f \wedge \Omega^1)_q,$$
will give again
$\dim E_{3}^{-1,2}(f)_{q _2}$ for $2d+1 \leq q \leq 3d $, resp. $\dim E_{3}^{-2,3}(f)_{q -3d}$ for $3d+1 \leq q \leq 4d $

\section{Computation of  $E^{s,t}_{\infty}(f)_k$ for $s+t=2$ under some assumptions}

We start with the following very useful result.

\begin{prop}
\label{propMS}
Conjecture \ref{conj30} holds if and only if 
$$\dim E_{3}^{2,0}(f)_k + \dim E_{3}^{1,1}(f)_k+\dim E_{3}^{0,2}(f)_k
-\dim E_{3}^{1,0}(f)_k - \dim E_{3}^{0,1}(f)_k+\delta_d^k=\chi(U),$$
for any integer $k \in [1,d]$, where $\delta_d^k=0$ if $k \ne d$ and $\delta_d^d=1$.
\end{prop}

\proof

Note that, using the properties of the filtrations $F$ and $P$ given in \eqref{PFincl} and \eqref{PFincl1}, it follows that $E^{s,t}_{\infty}(f)_k$, with $s+t=2$, can be nonzero only for $t=0,1,2$. Moreover, as explained above, one  has
$$\dim E_{3}^{1,0}(f)_k + \dim E_{3}^{0,1}(f)_k=\dim E_{\infty}^{1,0}(f)_k + \dim E_{\infty}^{0,1}(f)_k=\dim H^{1}(F,\C)_{\lambda },$$
where $\lambda= \exp(-2\pi i k/d)$.
On the other hand
$$\dim E_{3}^{2,0}(f)_k + \dim E_{3}^{1,1}(f)_k+\dim E_{3}^{0,2}(f)_k \geq $$
$$\dim E_{\infty}^{2,0}(f)_k + \dim E_{\infty}^{1,1}(f)_k+\dim E_{\infty}^{0,2}(f)_k =\dim H^{2}(F,\C)_{\lambda },$$
and the equality holds if and only if Conjectures \ref{conj30} holds. This equality is easy to check in practice as soon as we have the dimensions of the 3-terms  $E_{3}^{2,0}(f)_k,  E_{3}^{1,1}(f)_k$ and $E_{3}^{0,2}(f)_k$, since 
$$\dim H^{0}(F,\C)_{\lambda }-\dim H^{1}(F,\C)_{\lambda }+\dim H^{2}(F,\C)_{\lambda }=\chi(U)$$
for any $\lambda$, as follows from the formula \eqref{Euler}. Note that $\dim H^{0}(F,\C)_{\lambda }= \delta_d^k$, which completes our proof.

\endproof

\begin{rk}
\label{rkcheck}
Let $C: f=0$ be a reduced plane curve of degree $d$. To compute $\chi(U)$, it is enough to determine the Milnor number $\mu(C)$. This can be achieved in practice as follows. We choose the coordinates $x,y,z$ such that the line $z=0$ contains no singular point of $C$.
It follows that $\mu(C)=\mu(C_a)$, where $C_a$ is the affine curve in $\C^2$ defined by
$g(x,y)=f(x,y,1)=0$. Then one has 
$$\mu(C_a)= \dim \frac{\C[x,y]}{(g_x,g_y,g^2)}.$$
Here we have added $g^2$ to the Jacobian ideal $(g_x,g_y)$ since we want to take the sum of the Milnor numbers of the singularities of the polynomial $g$ situated on the curve $g=0$. For such a singularity $p$, the Brian\c con--Skoda Theorem, see \cite{BrSk}, shows that the local Milnor number at $p$ is not affected by the addition of the generator $g_p^2$, since 
$g_p^2 \in (g_{p,x}, g_{p,y})$. Here $g_p$ denotes the analytic germ induced by the polynomial $g$ at $p$.
\end{rk}

We assume mainly in this section that the curve $C:f=0$ has some non weighted homogeneous singularity, i.e. $\tau(f)=\tau(C) <\mu(C)=\mu(f)$. However, the formulas \eqref{vanishP} and \eqref{vanishP1} hold for any reduced plane curves satisfying the stronger Conjecture  \ref{conjmu-tau}, which explains why Corollaries \ref{corAA} and \ref{corAC} treat curves with weighted homogeneous singularities.
We also assume that
Conjecture \ref{conj30}  holds. When the stronger Conjecture  \ref{conjmu-tau} holds, the algorithm is much faster, since instead of computing the invariants till $q=4d$ we can stop practically at $q=q_0(f)\leq 3d+1$.

Note that for $t=0$ and $k=1,...,d$ we get the following
$$ \dim Gr_P^2 H^{2}(F,\C)_{\lambda }=\dim E^{2,0}_{\infty}(f)_k= \dim E^{2,0}_{3}(f)_k=$$
$$=\dim E^{2,0}_{2}(f)_k- \dim E^{0,1}_{2}(f)_k+\epsilon _{d+k} = \theta_k -\epsilon' _{d+k} +\epsilon _{d+k}.$$

When the stronger Conjecture  \ref{conjmu-tau} holds, and $t \geq 1$, one has also
$$ \dim Gr_P^{2-t} H^{2}(F,\C)_{\lambda }=\dim E^{2-t,t}_{\infty}(f)_k
=\dim E^{2-t,t}_{3}(f)_k=$$
$$=\dim E^{2-t,t}_{2}(f)_k-\mu(C)+\dim E_2^{-t,t+1}(f)_k=\theta_q -\epsilon'_{q+d},$$
for $t=1, 2$ and $k=1,2,..., d$, where we set $\epsilon'_q=\mu(C)-\tau(C)$ for $q>3d$. 
In particular, for  $q \geq q_1(f)=\max(2d-2, q_0(f))$, we have
$$\dim Gr_P^{2-t} H^{2}(F,\C)_{\lambda }=m(f)_{q-3}-\tau(C).$$
Note that, when $t=2$, this implies 
\begin{equation} 
\label{vanishP}
\dim Gr_P^{0} H^{2}(F,\C)_{\lambda }=m(f)_{2d+k-3}-\tau(C)=0
\end{equation} 
if $q \geq q_2(f)=\max(st(f)+3,q_1(f))$ and 
$$\dim Gr_P^{0} H^{2}(F,\C)_{\lambda }=m(f)_{2d+k-3}-\tau(C) \ne 0$$
if $q = st(f)+2 \geq q_1(f)$.
Similarly, for $t=1$ and $k=d-2$ or $k= d-1$ we get
\begin{equation} 
\label{vanishP1}
\dim Gr_P^{1} H^{2}(F,\C)_{\lambda }=m(f)_{d+k-3}-\tau(C)=0,
\end{equation} 
if $q =d+k\geq q_2(f).$
\begin{cor}
\label{corAA}
Let $\A: f=0$ be an arrangement of $d$ lines in $\PP^2$ with Milnor fiber $F$.
Then $P^1H^2(F,\C)=H^2(F,\C)$ and $Gr_P^{1} H^{2}(F,\C)_{\lambda }=0$ for
$\lambda=\exp(-2\pi i (d-1)/d)$.
\end{cor}

\proof
The first part of the claim is clearly equivalent to $\dim Gr_P^{0} H^{2}(F,\C)_{\lambda }=0$ for any $\lambda$.
Since $st(f) \leq 2d-4$ for any line arrangement, see \cite{DIM2}, the claim follows from 
\eqref{vanishP}. The second part of the claim follows from \eqref{vanishP1}.
\endproof

\begin{cor}
\label{corAB}
Let $C: f=0$ be a reduced plane curve which is either free or nearly free, and for which Conjecture \ref{conjmu-tau} hold, with $q_0(f) \leq 2d+1$.
Then $P^1H^2(F,\C)=H^2(F,\C).$
\end{cor}

\proof
When  $C:f=0$ is a free (resp. nearly free) curve with exponents $d_1 \leq d_2$, it is known that $st(f)=d_2+d-3 \leq 2d-4$ (resp. $st(f)=d_2+d-2 \leq 2d-3$. We conclude as above.
\endproof

\begin{cor}
\label{corAC}
Let $C: f=0$ be a reduced plane curve of degree $d\geq 5$ having only one singularity, which is a node.
Then $\Delta_C^1(t)=1$ and
$$\dim Gr_P^{0} H^{2}(F,\C)_{\lambda }=m(f)_{2d+k-3}-\tau(C)=m(f)_{d-k-3}-1={d-k-1 \choose 2}-1 \ne 0$$
for $k=1,2,...,d-4.$
\end{cor}

\proof For the first claim, see for instance \cite[Theorem 6.4.17]{D1}. To prove the second claim,
note that for such a curve $ct(f)=st(f)=3d-6$, see \cite{DStEdin}. This implies in particular
$$m(f)_{2d+k-3}=m(f_s)_{2d+k-3}=m(f_s)_{d-k-3}=m(f)_{d-k-3},$$
where the equality in the middle follows from the Grothendieck duality of the Artinian  Milnor algebra $M(f_s)$.
\endproof
Note that for this uninodal curve, the vanishing bound given in \eqref{vanishP} is sharp.
For more on such uninodal curves, see Corollary \ref{corAC'} below.

\section{Application to the study of Bernstein-Sato polynomials}
Let $(D,0):g=0$ be a complex analytic hypersurface germ at the origin of $\C^n$ and denote by
$b_{g,0}(s)$ the corresponding (local) Bernstein-Sato polynomial. If the analytic germ $g$ is given by a homogeneous polynomial, then one can define also the global Bernstein-Sato polynomial $b_g(s)$ of $g$,
and one has an equality $b_g(s)=b_{g,0}(s)$, see for more details \cite{Sa1}, \cite{Sa2}.
Let $R_{g,0}$ (resp. $\tilde R_{g,0}$) be the set of roots of the polynomial $b_{g,0}(-s)$ (resp. of the polynomial $b_{g,0}(-s)/(1-s)$). Recall that one has
$$\tilde R_{g,0} \subset [\tilde \alpha_{g,0}, n-\tilde\alpha_{g,0}],$$
where $\tilde\alpha_{g,0} =\min \tilde R_{g,0}>0$. Moreover, $\alpha_{g,0} =\min R_{g,0}>0$ coincides with the log canonical threshold $lct(g)$ of the germ $g$, a.k.a. the log canonical threshold $lct(\C^n,D)$ of the pair $(\C^n,D)$, see \cite{Ko}.
When $g$ is a homogeneous polynomial, we use the simpler notation $R_g=R_{g,0}$, $\tilde R_g=\tilde R_{g,0}$ and so on.

\begin{ex}[Cones over smooth projective hypersurfaces]
\label{exBS}
Let $g$ be a homogeneous polynomial of degree $d$ in $n$ variables such that the corresponding projective hypersurface $g=0$ is smooth.
Using the relation between the zero set $\tilde R_g$ and the spectrum in the case of an isolated weighted homogeneous hypersurface singularity, see for instance \cite{SaSurvey}, one has 
$$\tilde R_g=\left \{\frac{j}{d} \ : \ n \leq j \leq n(d-1)\right \}.$$

\end{ex}

In this section we consider mainly the case when $n=3$ and $g=f$ is the defining equation of a reduced curve $C$ in $\PP^2$. Let $D$ be the surface in $\C^3$ defined by $f=0$ and note that at a point $a \in D$, $a \ne 0$, the germ $(D,a)$ is analytically a product between a plane curve singularity and a smooth 1-dimensional germ. It follows that 
$$\tilde R_{f,a} \subset [\tilde \alpha_{f,a}, 2-\tilde \alpha_{f,a}].$$
Recall  M. Saito's  fundamental result \cite[Theorem 2]{Sa1} quoted here in the case $n=3$.

\begin{thm}
\label{thmBS}
Let $C: f=0$ be a reduced curve in $\PP^2$, let $\alpha >0$ be a rational number 
and set $\lambda=\exp(-2 \pi i \alpha)$. 

\begin{enumerate}

\item If $Gr_P^p H^{2}(F,\C)_{\lambda }\ne 0$, where $p=\lfloor 3-\al \rfloor $, then $ \al \in R_f$

\item If the sets $\al+\N$ and  $ \cup_{a \in D, a \ne 0}R_{f,a}$ are disjoints,
then the converse of the assertion $(1)$ holds.
\end{enumerate}
\end{thm}
The following result is due to M. Saito in arbitrary dimension, see \cite[Theorem 1]{Sa0}.
We give a new, simple proof below to point out the relation with Corollary \ref{corAA}.
\begin{cor}
\label{corAA'}
Let $\A: f=0$ be an arrangement of $d$ lines in $\PP^2$ with Milnor fiber $F$.
Then 
$$\max R_f <2-\frac{1}{d}.$$
\end{cor}
\proof We prove the result by discussing the possible cases for $\al$.
If $\al=2$, the result is clear by Theorem \ref{thmBS}, since 
$$H^{2}(F,\C)_1=H^2(U, \C)=F^2H^2(U, \C) = F^2H^{2}(F,\C)_1 \subset P^2H^{2}(F,\C)_1 $$
and hence $H^{2}(F,\C)_1=P^2H^{2}(F,\C)_1 $, which implies $Gr_P^1 H^{2}(F,\C)_{1}= 0$.
If $\al \ne 2$, then we know that $d \al \in \N$, and hence to prove the claim we have to consider
the case $\al'=1+ (d-1)/d$ and the case $\al''=2+k/d$ for $k=1,..., d-1$.
Both cases follow from Theorem \ref{thmBS} (2) and Corollary \ref{corAA}. To see this, the only point to explain is why
$$\al' \notin  \cup_{a \in D, a \ne 0}R_{f,a}.$$
Indeed, note that a singular point $a \in D$ correspond to a point of multiplicity $m \geq 2$ in the line arrangement. For such a point one has $\alpha_{f,a}=2/m$, and this implies our claim.
\endproof
The following result can be proved along the same lines using Corollary \ref{corAB}. \begin{cor}
\label{corAB'}
Let $C: f=0$ be a reduced plane curve which is either free or nearly free, and for which Conjectures \ref{conj3} and  \ref{conjmu-tau} hold.
Then 
$\max R_f \leq 2.$
\end{cor}
The next result shows that a nodal curve behaves quite differently from a line arrangement with respect to the zero set $R_f$.
\begin{cor}
\label{corAC'}
Let $C: f=0$ be a reduced plane curve of degree $d\geq 5$ having only one singularity, which is a node.
Then 
$$ R_f=\left \{\frac{j}{d} \ : \ 3 \leq j \leq 3d-4\right \}.$$
In particular
$$\max R_f = 2+ \frac{d-4}{d}.$$
\end{cor}

\proof Recall that for a uninodal plane curve $C:f=0$ of degree $d$ one has the following: $ct(f)=st(f)=3d-6$ and $\dim H^2(K_f^*)_q=0$ for $q<2d-2$ and $\dim H^2(K_f^*)_q=1$ for $q \geq 2d-2$, see \cite[Example 4.3 and Equation (2.17)]{DStEdin}. Using this, the fact that $R_f$ contains the given set of rational numbers follows Theorem \ref{thmconj} and Theorem \ref{thmBS} (1). To show that $\al=(3d-3)/d$ is not in $R_f$ one uses 
Theorem \ref{thmBS} (2).

\endproof

\begin{rk}
\label{rkBS}
\begin{enumerate}

\item  Corollary \ref{corAC'} can be restated as
$$R_f=R_{f_s} \setminus \left \{\frac{d-3}{d}\right \},$$
where $f_s$ is a generic polynomial of degree $d$ in $x,y,z$ as in Example \ref{exBS}.

\item In the case of a line arrangement $\A:f=0$ in $\PP^2$, the zero set $R_f$ is not determined by the combinatorics, see Walther \cite{Wa} and Saito \cite{Sa0}. In fact, there is a pair of line arrangements $\A_1:f_1=0$ and $\A_2:f_2=0$ of degree $d=9$, going back to Ziegler \cite{Zieg}, having the same combinatorics but different sets $R_f$ and different Hilbert functions for their Milnor algebras $M(f_1)$ and $M(f_2)$.

\item The zero set $R_f$ is not determined by the Hilbert function for its Milnor algebra $M(f)$ for a reduced plane curve $C:f=0$, see \cite{Sa2} and Example \ref{exSaito} below.

\end{enumerate}
\end{rk}

\section{Examples}

\begin{ex}[A torus curve of type $(2,4)$]
\label{extc2,4}

Consider the curve $C: f=(x^2+y^2)^4+(y^4+z^4)^2=0$. This curve has 8 weighted homogeneous singularities of type $A_3$ with local equation $u^2+v^4=0$, and hence
$\mu(C)=\tau(C)=8 \times 3=24.$
A direct computation shows that $\epsilon '_q=0$ for $q\geq q_0(f)=11$ as in Conjecture \ref{conjmu-tau}. On the other hand we have
$\theta_q=0$ for $q \geq 18$, which is exactly the bound predicted by \eqref{vanishP} which is in this case $q \geq q_2(f)=st(f)+3=15+3=18.$
To state the full result, we consider the pole order spectrum defined by
\begin{equation} 
\label{sp1}
Sp_P^j(f)=\sum_{\al>0}n^j_{P,f,\al}t^{\al} 
\end{equation} 
for $j=0,1$, where 
$$n^j_{P,f,\al}=\dim Gr_P^pH^{2-j}(F,\C)_{\lambda}$$
with $p=[3-\al]$ and $\lambda=\exp(-2 \pi i \alpha)$. With this notation, we have the following for our torus curve of type $(2,4)$.
\begin{equation} 
\label{sp11}
Sp_P^1(f)=t^{\frac{14}{8}}+t^{\frac{16}{8}}+t^{\frac{18}{8}}.
\end{equation} 

\begin{equation} 
\label{sp10}
Sp_P^0(f)=t^{\frac{3}{8}}+3t^{\frac{4}{8}}+6t^{\frac{5}{8}}+10t^{\frac{6}{8}}+12t^{\frac{7}{8}}+16t^{\frac{8}{8}}+18t^{\frac{9}{8}}+20t^{\frac{10}{8}}+
\end{equation} 
$$+18t^{\frac{11}{8}}+16t^{\frac{12}{8}}+13t^{\frac{13}{8}}+10t^{\frac{14}{8}}+7t^{\frac{15}{8}}+3t^{\frac{16}{8}}+t^{\frac{17}{8}}.$$
\end{ex}
Note in particular that, unlike the case of line arrangements treated in Corollary \ref{corAA},
here $P^1H^2(F,\C) \ne H^2(F,\C)$.

\begin{ex}[A free curve with non weighted homogeneous singularities]
\label{exfreem=2}
Consider the curve $C: f=(y^2z^2-x^{4})^2y^2-x^{10},$
 introduced in \cite{B+}. This curve is free with exponents $(4,5)$, and hence $st(f)=12.$
Moreover, one has $\mu(C)=70$ and $\tau(C)=61$. A direct computation shows that $\epsilon '_q=\mu(C)-\tau(C)=9$ for $18= q_0(f) \leq q \leq 30$. Moreover, we have
$ Gr_P^{1} H^{2}(F,\C)_{\lambda }=0$ for $q =d+k \geq 15$, which is stronger that the bound predicted by \eqref{vanishP1} which is in this case $q=d+k \geq d+(d-2)=18 \geq st(f)+3=15.$ Moreover Conjecture \ref{conj30}  holds for this curve since the condition in Proposition \ref{propMS} is fulfilled. 
In this case we have
\begin{equation} 
\label{sp21}
Sp_P^1(f)=t^{\frac{16}{10}}+t^{\frac{17}{10}}+t^{\frac{18}{10}}  +t^{\frac{19}{10}} +t^{\frac{20}{10}} +t^{\frac{21}{10}} +t^{\frac{22}{10}} +t^{\frac{23}{10}} +t^{\frac{24}{10}} 
\end{equation} 
and
\begin{equation} 
\label{sp22}
Sp_P^0(f)=t^{\frac{3}{10}}+2t^{\frac{4}{10}}+3t^{\frac{5}{10}}+4t^{\frac{6}{10}}+4t^{\frac{7}{10}}+4t^{\frac{8}{10}}  +4t^{\frac{9}{10}} +3t^{\frac{10}{10}} +4t^{\frac{11}{10}} +4t^{\frac{12}{10}} +3t^{\frac{13}{10}} +2t^{\frac{14}{10}} .
\end{equation} 
Moreover, $P^1H^2(F,\C)=H^2(F,\C)$ as predicted by Corollary \ref{corAB}.

The free curve $C$ above has two irreducible components. To get a similar example with an irreducible free curve, one may consider 
$$C':f'=x^4y^2+y^6-3xy^4z+3x^2y^2z^2 -x^3z^3=0.$$
Then $C'$ is free with exponents $(2,3)$ and one has $\mu(C')=20$, $\tau(C')=19$, $q_0(f)=10$. Our algorithm gives
\begin{equation} 
\label{spf'1}
Sp_P^1(f')=t^{\frac{11}{6}}+t^{\frac{13}{6}}
\end{equation}  
and, respectively,
\begin{equation} 
\label{spf'0}
Sp_P^0(f')=t^{\frac{3}{6}}+t^{\frac{4}{6}}+2t^{\frac{5}{6}}+2t^{\frac{7}{6}}+t^{\frac{8}{6}}.
\end{equation}

\end{ex}

\begin{ex}[A curve with  $\mu(C) \ne \tau(C)$ and high $q_0(f)=2d+2$]
\label{exhighq0}

Consider the \\
curve $C: f=x^4y^4z^4+x^{12}+ y^{12}=0$.
 This curve is far from being free, since $ct(f)+st(f)-T=11$, where $T=3(d-2)$, and it has $\mu(C)=73$ and $\tau(C)=64$. 
 A direct computation shows that $\epsilon '_q=\mu(C)-\tau(C)=9$ for $q\geq q_0(f)=26=2d+2<3d$, as required in Conjecture \ref{conjmu-tau}. In fact, Conjecture \ref{conj30}  holds for this curve since the condition in Proposition \ref{propMS} is fulfilled. 
 Note also that in this example the sequence 
$\epsilon '_q$ contains the sequence  $\epsilon '_{20}=3$, $\epsilon '_{21}=\epsilon '_{22}=2$,
$\epsilon '_{23}=3$, and hence it is not, as in most other examples, an increasing sequence.
 Note that $H^1(F,\C)=0$ in this case,  $b_2(F)=455$ is large, and $P^1H^2(F,\C)=H^2(F,\C)$.

\end{ex}

\begin{ex}[Some examples considered by M. Saito in \cite{Sa2}]
\label{exSaito}
Consider the curves $C_1: f_1=x^5+y^4z+x^4y=0$ and $C_2:f_2=x^5+y^4z+x^3y^2=0$.
Then the Milnor algebras $M(f_1)$ and $M(f_2)$ have the same Hilbert functions, and one has
$\mu(f_1)=\mu(f_2)=12$ and $\tau(f_1)=\tau(f_2)=11$. Our algorithm shows that in  both cases $H^1(F,\C)=0$, $q_0(f_1)=9$, $q_0(f_2)=8$, while the spectrum on $H^2(F,\C)$ is given by
\begin{equation} 
\label{spf1}
Sp^0(f_1)=t^{\frac{3}{5}}+t^{\frac{4}{5}}+t^{\frac{6}{5}}+t^{\frac{7}{5}},
\end{equation} 
and, respectively,
\begin{equation} 
\label{spf2}
Sp^0(f_2)=t^{\frac{4}{5}}+t^{\frac{6}{5}}+t^{\frac{7}{5}}+t^{\frac{8}{5}},
\end{equation}  
hence exactly the same formulas as in \cite{Sa2}, equations (1.4.5). Using Theorem \ref{thmBS}, it follows that $$ \frac{8}{5} \in R_{f_2} \setminus R_{f_1}$$
as proved in \cite{Sa2}. Notice also that $C_1:f_1=0$ and $C_2:f_2=0$ are nearly free with exponents $(2,3)$. Both curves have a unique unibranch singulary located at $p=(0:0:1)$, which is semi-weighted-homogeneous with weights $(1/5,1/4)$. In particular, the singularities $(C_1,p)$ and $(C_2,p)$ are topologically equivalent.

In \cite[Remark 2.5]{Sa2}, M. Saito also considers the curves
$C_3:f_3=x^5+xy^3z+y^4z+xy^4=0$ and $C_4:f_4=x^5+xy^3z+y^4z=0$. 
In this case one has $\mu(f_3)=\mu(f_4)=11$,  $\tau(f_3)=\tau(f_4)=10$. Our algorithm shows that again $H^1(F,\C)=0$, $q_0(f_3)=10$, $q_0(f_4)=9$, while the spectrum on $H^2(F,\C)$ is given by
\begin{equation} 
\label{spf3}
Sp^0(f_3)=t^{\frac{3}{5}}+2t^{\frac{4}{5}}+t+2t^{\frac{6}{5}}+2t^{\frac{7}{5}}+t^{\frac{8}{5}},
\end{equation} 
and, respectively,
\begin{equation} 
\label{spf4}
Sp^0(f_4)=t^{\frac{3}{5}}+ t^{\frac{4}{5}}+t+ 2t^{\frac{6}{5}}+2t^{\frac{7}{5}}+t^{\frac{8}{5}}+t^{\frac{9}{5}}.
\end{equation}  
Using Theorem \ref{thmBS}, it follows that $$ \frac{9}{5} \in R_{f_4} \setminus R_{f_3}$$
as proved in \cite{Sa2}.
Notice also that $C_3:f_3=0$ and $C_4:f_4=0$ are neither free nor nearly free, they satisfy
$ct(f_j)+st(f_j)-T=3$ for $j=3,4$. In both cases, there are independent syzygies of degree 2 and 3. Note that Conjecture \ref{conj30}  holds for all these curves, since the condition in Proposition \ref{propMS} is fulfilled. Both curves have a unique two-branch singulary located at $p=(0:0:1)$, which is semi-weighted-homogeneous with weights $(1/5,4/15)$. Again, the singularities $(C_3,p)$ and $(C_4,p)$ are topologically equivalent.

\end{ex}

\end{document}